\documentclass{article}
\usepackage[utf8]{inputenc}
\usepackage[a4paper,top=2.5cm,bottom=2.5cm,left=2.5cm,right=2.5cm,marginparwidth=1.75cm]{geometry}
\usepackage{amssymb}
\usepackage{amsmath}
\usepackage{amsthm}
\usepackage{epsfig}
\usepackage{graphicx}
\usepackage{graphics}
\usepackage{color}
\usepackage{authblk}
\usepackage{hyperref}
\usepackage{pdfpages}
\usepackage{soul}

\begin{document}

\title{Dynamic Programming Techniques \\ for Planar Orbital Transfer of Low Earth Orbit Satellites}
\author[1]{C. Ciancarelli}
\author[2]{R. Ferretti}
\author[1]{A. Intelisano}
\author[3]{G. Villani}

\affil[1]{Thales Alenia Space Italia, Roma (Italy)}
\affil[2]{Dipartimento di Matematica e Fisica, Universit\`a Roma Tre, Roma (Italy)}
\affil[3]{Dipartimento di Matematica, ``Sapienza'' Universit\`a di Roma and Thales Alenia Space Italia, Roma (Italy)}

\date{\today}

\maketitle

\begin{abstract}

In this paper, we present an application of the optimal control theory to orbital transfer of Low Earth Orbit satellites. The optimal control problem is treated with Dynamic Programming techniques which require solving the Hamilton--Jacobi--Bellman equation on a suitable state space, with the reconstruction of the optimal controls made in the form of a static feedback. In order to validate the numerical scheme without the complexity of the full model, this first study sets the problem in planar form, thus working in a four-dimensional state space. We will study various techniques to speed up the computation, and assess the accuracy of the numerical solution.

This project is born from the attempt of evaluating and applying direct method of optimal control techniques based on Dynamic Programming as a complementary approach to the well known indirect methods, as Pontryagin or Lawden. In particular, the final aim is to treat the case of low thrust engines from real use cases, in the full 3-D problem.

\end{abstract}



\section{Introduction}

A Low Earth Orbit (LEO) is an orbit that is relatively close to Earth’s surface. It usually lies at an altitude of less than 1000 km but could be as low as 160 km above Earth -- which is low if compared to other orbits, but still very far above Earth’s surface.
Unlike satellites in Geostationary Earth Orbit (GEO) that must always orbit along Earth’s equator, LEO satellites might have a tilted orbit plane. This means there are more available routes for satellites in LEO, which is one of the reasons for LEO to be a very commonly used orbit.

LEO’s close proximity to Earth makes it useful for several reasons. It is the orbit most commonly used for satellite imaging, as being near the surface it allows the satellite to take images at a higher resolution. It is also the orbit used for the International Space Station (ISS), being it easier for astronauts to travel to and from it due to the shorter distance. Satellites in this orbit travel at a speed of about 7.8 km per second; at this speed, a satellite takes approximately 90 minutes to circle Earth.


%
%
%

Our goal here is to optimize the orbit raising manoeuvre of the orbital transfer between the preliminary orbit and the operative orbit, focusing in particular on small LEO satellites with low thrust, electric propulsion. This target will be pursued by Dynamic Programming (DP) techniques, keeping in mind that this approach might eventually be integrated with more usual techniques like  Pontryagin Maximum Principle (PMP).

The use of DP techniques in space systems mission analysis is not common, due to their high requirements in terms of computer resources. To our knowledge, a line of research in this direction has been pursued by Bokanowski et al. in a series of papers (see \cite{BBDZ15,BBDZ21,BBZD} and the references therein), devoted to the global optimization of the ascent trajectory of a launcher for geostationary satellites. Although our work shares some techniques with theirs, we are interested here in obtaining an optimal solution in the form of a static feedback, as well as in cutting the complexity of numerical computations in order to avoid the use of massively parallel machines.

Dynamic Programming and PMP have complementary features. DP allows for the computation of approximate optimal controls in feedback form, while PMP provides optimal controls in open-loop form; DP gives a global minimum, while PMP may provide a local one. On the other hand, DP suffers from the ``curse of dimensionality", and for this reason, if still feasible, it typically produces solutions which are expected to be less accurate than the corresponding solutions obtained via the PMP. Thus, one typical way of complementing DP and PMP is to use the less accurate (but global) DP solution as a first guess to initialize PMP and avoid as much as possible convergence to local minima. We refer to \cite{CriMar} for an extensive study of the integration between these two techniques; a comparison between the two approximations of the optimal solutions is also carried out in \cite{BBDZ15,BBDZ21,BBZD} for the case of the launcher trajectory.

In this paper, we study techniques to minimize the complexity issue of the DP technique, working on a simplified 2-D version of the specific problem under consideration. In particular, this work integrates:

\begin{itemize}

\item
A reduced space grid to save memory space;

\item
A Policy Iteration (PI) algorithm to accelerate convergence for the iterative computation of the value function;

\item
A fast direct search algorithm to perform minimization.

\end{itemize}

In addition, we study a proper definition of the cost functional, which could allow to obtain the correct behaviour of the optimal solution.

The paper is structured as follows. Section \ref{sec:model} states the basic model for the orbit transfer problem; Sect. \ref{sec:implementation} discusses the implementation details of the scheme and Sect. \ref{sec:test} shows some numerical tests on realistic cases.

\section{Mathematical Model}\label{sec:model}

We recall in this section the general setup for the control problem. Dynamic Programming leads to characterize the optimal solution through a Hamilton--Jacobi--Bellman (HJB) equation. In our case, we are not given a time horizon for the orbit transfer; in addition, we aim at obtaining a static feedback. These features are naturally implemented with an infinite horizon control problem, which results in a stationary HJB equation. We start by reviewing the dynamics of the satellite system.

\subsection{Dynamics of a satellite}
For the formulation of the mathematical model of a satellite we consider the total acceleration:
\begin{equation}\label{eq:model0}
    \ddot{\Vec{r}} = - \frac{\mu}{r^3}\Vec{r} + \Vec{a}_{p} + \Vec{a}_{T}.
\end{equation}
In the right-hand side, the first term is the gravity acceleration, while $\Vec{a}_{T}$ represents the acceleration generated by the satellite thrusters and $\Vec{a}_{p}$ is the total acceleration due to the perturbations effect, and more precisely
\begin{equation*}
    \Vec{a}_{p} = \Vec{a}_{ns} + \Vec{a}_{sun} + \Vec{a}_{moon} + \Vec{a}_{D} + \Vec{a}_{SR} + \Vec{a}_{E}
\end{equation*}
\begin{itemize}
    \item  $\Vec{a}_{ns}$ due to the asphericity of the Earth
    \item  $\Vec{a}_{sun}$ and $\Vec{a}_{moon}$ for solar and lunar attraction respectively
    \item  $\Vec{a}_{D}$ due to aerodynamic drag
    \item  $\Vec{a}_{SR}$ for the solar radiations
    \item  $\Vec{a}_{E}$ due to the electromagnetic effects.
\end{itemize} 
The effect of each component of $\Vec{a}_{p}$ on the dynamics of the satellite depends on the altitude of the satellite as well as on time. 
In order to formulate the problem in terms of a stationary HJB equation, all perturbations appearing in the term $\Vec{a}_{p}$, and depending on time, will not be taken into consideration. We will eventually introduce an approximate aerodynamic drag in the numerical tests, to show that the optimal feedback generated by the scheme under consideration tends to correct the effects of omitted terms.

As for the control term $\Vec{a}_{T}$, we will consider the case in which the thrusters can only work at full thrust, or being switched off. If working at full thrust,  $\Vec{a}_{T}$ can take any direction, and, in order to avoid the introduction of further state variables (which would drastically increase the complexity of the problem), we will neglect the effects of rotational inertia and assume that the orientation of the thruster is instantaneous. The variation of mass due to fuel consumption will be neglected as well.

In the simplified $\mathbb{R}^2$ version, the dynamics is recast more explicitly as
the following controlled system:
\begin{equation}\label{eq:R2cartesiano}
    \begin{cases}
        \dot{y_1} = v_1 \\
        \displaystyle \dot{v}_1 = - \frac{\mu y_1}{(y_1^2+y_2^2)^{3/2}}  + u_1 \\
        \dot{y_2} = v_2 \\
        \displaystyle \dot{v}_2 = - \frac{\mu  y_2}{(y_1^2+y_2^2)^{3/2}} + u_2,
    \end{cases}
\end{equation}
with the initial conditions
\begin{equation*}
     \begin{cases}
        y_1(0) = {y_1}_0 \\
        v_1(0) = {v_1}_0 \\
        y_2(0) = {y_2}_0 \\
        v_2(0) = {v_2}_0,
     \end{cases}
\end{equation*}
where $\mu$ is the gravitational parameter associated to the main attracting body,  for the Earth $\mu = 398600.4 $ km$^3$/s$^2$, and the controls $u_1$ and $u_2$ are the two components of the acceleration provided by the satellite thrusters. For any given time $t>0$, the control $u(t)=(u_1(t),u_2(t))$ belongs to the set
\begin{equation*}
U = \{u_1,u_2: u_1^2+u_2^2=T^2/m^2\} \cup \{0\},
\end{equation*}
where $T$ is the maximum thrust and $m$ is the satellite mass.

\subsection{Infinite horizon problem}
Consider now a control system governed by the state equation
\begin{equation}\label{eq:dyn}
    \begin{cases}
        \dot y(x,t,u) = b(y(t), u(t)), \ \ t>0, \\
        y(0) = x.
    \end{cases}
\end{equation}
Here, the control $u$ is a measurable function of $t \in [0,+\infty)$ with values in a closed bounded set in $\mathbb{R}^M$, and $b(\cdot,\cdot)$ is assumed to be Lipschitz continuous w.r.t. both variables. In the model under consideration, $y=(y_1,v_1,y_2,v_2)$, $x=({y_1}_0,{v_1}_0,{y_2}_0,{v_2}_0)$, $u=(u_1,u_2)$, and the set of admissible control functions is
\begin{equation}
\mathcal{U} = \{u:[0,+\infty)\to U, u \text{ measurable}\}
\end{equation}
The cost functional to be minimized has the general form
\begin{equation}\label{eq:cost}
    J(x,u) = \int_{0}^{+\infty} l(y(x,t,u),u(t))e^{-\lambda t} dt,
\end{equation}
where $\lambda > 0$ represents a fixed discount factor, and $l$ is the so-called {\em running cost}. The first step in the standard DP approach to the optimal control problem described above is to introduce the value function $V$, defined by
\begin{equation}\label{eq:valuefunction}
    V(x) := \inf_{u \in \mathcal{U}} J(x,u).
\end{equation}
Then, it can be proved that $V$ is a weak solution, in the viscosity sense (see \cite{BCD,FF13}), of the HJB equation
\begin{equation}\label{eq:HJB}
    \lambda V(x) + \sup_{f \in U}\{-b(x,f)\cdot DV(x) - l(x,f) \} = 0.
\end{equation}
Some remarks should be made about choosing the cost functional \eqref{eq:cost}. First, we are treating orbit transfer problems with no fixed time horizon -- this is reasonable, in particular, for satellites with low-thrust engines. Second, the form \eqref{eq:dyn}--\eqref{eq:cost} gives an optimal solution which can be expressed in terms of a static feedback, i.e., at any given time $t$ the optimal control $u^*(t)$ is constructed as
\begin{equation}\label{eq:feedback}
u^*(t) = F(y^*(t)),
\end{equation}
where $F$ is the feedback function (in the discrete setting, we will also term this function as {\em policy}), and $y^*$ is the trajectory associated to $u^*$. As far as the supremum in \eqref{eq:HJB} is a maximum, the value of $f$ achieving this maximum should be regarded as the value $F(x)$ of the feedback function at $x$.

The possibility of constructing the optimal solution by a static feedback may be useful in case of any unexpected perturbation of the optimal trajectory, which can be corrected without computing a new optimal solution from scratch. In order to write in problem in the form \eqref{eq:dyn}--\eqref{eq:cost}, we are forced to assume a dynamics in the time-independent form \eqref{eq:dyn}, and therefore to neglect all perturbations depending on time, as remarked above.

\paragraph{Choice of the running cost}

To reach the operative orbit we should properly define the targets to be achieved. In our case, the running cost must take into account the control, i.e., the thrust applied from the engines of the satellite, and the orbital parameters which uniquely identify the target orbit; when used in the cost functions, these parameters should be discriminating enough to force a good agreement between the final and the target orbit. We found that a suitable set of parameters for the running cost is given by eccentricity $e$,  semi-major axis $a$,  argument of periapsis $\omega$ (in this case, argument of perigee), and true anomaly $\theta_*$ (in three dimensions, we should also take into account the inclination $i$ and longitude of the ascending node $\Omega$).
We consider therefore a running cost in the following compact notation
\begin{equation}\label{eq:runcost}
    l(y,v,u) = \alpha |u| + \beta (a - \bar{a})^2 + \gamma \left((e \cos(\omega) - \bar{e}\cos(\bar{\omega}))^2 + (e \sin(\omega) - \bar{e}\sin(\bar{\omega})^2\right),
\end{equation}
where $a$ and $e$ are the semi-major axis and the (scalar) eccentricity of the orbit corresponding to the state $(y,v)$, $\bar{a}$ and $\bar{e}$ are the semi-major axis and the eccentricity for the target orbit, and $|u|=\sqrt{u_1^2+u_2^2}$ is the norm of the control.
A control term in this form corresponds to what in classical optimal control is termed as {\em minimum fuel problem}; in our case, we look for a compromise between low fuel consumption and precise attainment of the target orbit. Note that, in \eqref{eq:runcost}, the eccentricity has been split in the two cartesian components $e_1=e\cos\omega$ and $e_2=e\sin\omega$, which also depend on $\omega$.
Since the scheme will work in polar coordinates, we will defer the explicit expression of the orbit parameters and write them directly in polar form. Note that an optimal control based on orbital parameters (in this case, to define the target set) is also studied in \cite{BBZD}.


A suitable, delicate balance of the various terms in \eqref{eq:runcost} should be pursued in order to obtain the required behaviour from the satellite system. In particular, if $\beta$ is too small, the optimal solution tends to be more precise in achieving the target orbit at the expense of fuel consumption (the thruster keeps on correcting the orbit even if the satellite is in its neighbourhood), while it shows an opposite behaviour if this parameter is too large (the thrusters eventually stop, but the operative orbit is achieved with a low accuracy).

\section{Implementation of the Numerical Scheme}\label{sec:implementation}

The HJB equation has been discretized in semi-Lagrangian (SL) form (see \cite{FF13}). Once a space grid with nodes $x_j$ has been set, this approach provides a scheme in the fixed point form
\begin{equation}\label{eq:schema}
v_j = \min_{f\in U} \left\{\Delta t \> l(x_j,f) + e^{-\lambda\Delta t} I[V](x_j+\Delta t\>\Phi(x_j,f))\right\} \qquad (j=1,\ldots,N)
\end{equation}
Here, $\Delta t$ is a (pseudo-)time discretization step, $x_j+\Delta t\>\Phi(x_j,f)$ is a one-step approximation of the trajectory starting at $x_j$ and subject to a control $f$, and $I[V](z)$ is a space interpolation of the values in $V$, computed at the point $z$. In this work, the approximation of trajectory has been performed with a second-order Runge--Kutta scheme, and the interpolation has been chosen in the multilinear ($\mathbb Q_1$) form. It is well-known (see \cite{FF13}) that this latter choice leads to a monotone scheme, and therefore, via the Barles--Souganidis theory \cite{BS}, ensures convergence to the viscosity solution of \eqref{eq:HJB}.

As customary in DP problems, any practical implementation must face a high complexity, in terms of both CPU time and memory requirements. The main bottlenecks of the scheme \eqref{eq:schema} are memory space required, efficiency of the iterative solution and minimization phase. We will examine these three issues in detail.

\begin{figure}[t!]
    \centering
    \includegraphics[width=9cm]{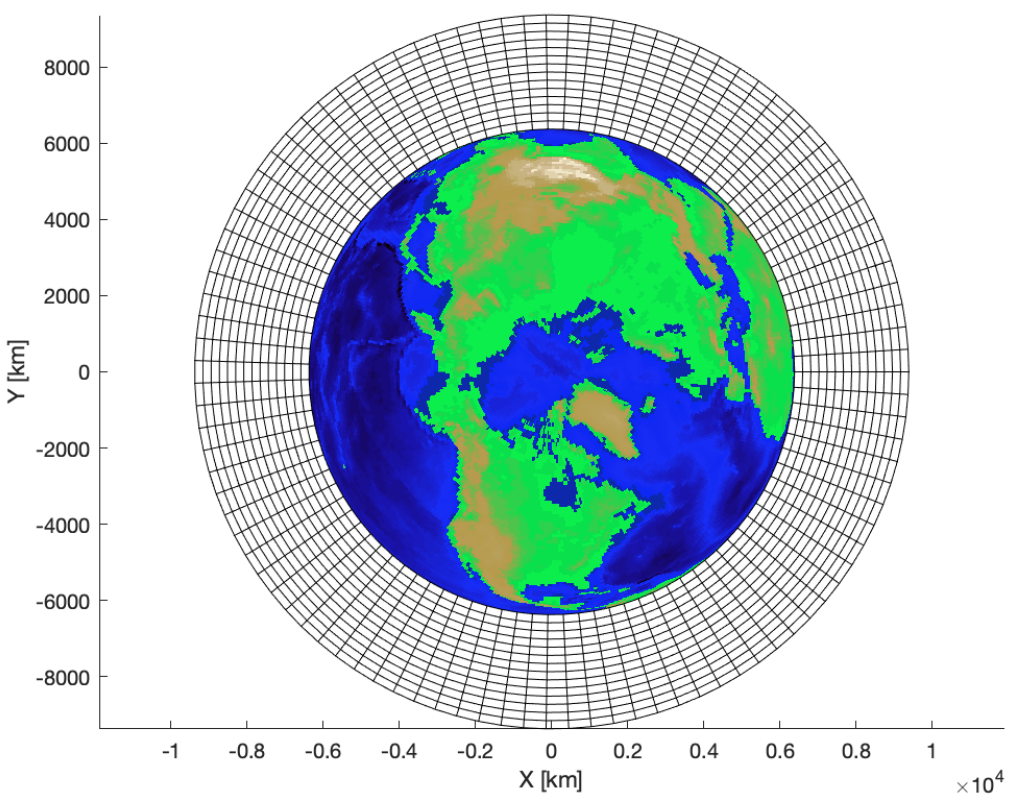}
    \caption{Circular crown grid}
    \label{fig:circular_crown}
\end{figure}

\paragraph{Grid and polar coordinates.}
The first trick to reduce memory requirements and number of mesh nodes is the choice of a grid in form of a circular crown, as in Fig. \ref{fig:circular_crown}, constructed by ``bending'' a rectangular grid, and closing it with periodic conditions; this is done not only in the space, but also in the velocity coordinates. In this way one avoids to compute the value function at grid points that may be far from the target, or even under the Earth surface\footnote{Note that, for the system under consideration, this construction avoids the origin, where a blow up of the dynamics occurs, and restricts the state space to a region in which the dynamics is actually Lipschitz continuous.}. Since we expect that optimal trajectories remain in a relatively close neighbourhood of the target orbit, the two circular portions of the boundary are provided with a state constraint condition, obtained by penalizing trajectories that fall outside of the domain, for examples by setting $I[V](z)=M$ (a sufficiently large value) as a stopping cost for $z$ out of the circular crown. This reduction of the computational domain, except for a somewhat different geometry, is much in the same spirit of \cite{BBDZ21}; in fact, this is a key step to reduce the computational resources required, and hence improve the applicability of DP techniques in an industrial framework.

Note that, with this construction, the grid results from the composition of a rectangular grid with a smooth transformation, and therefore the interpolation enjoys the same efficiency of the pure rectangular case.

In order to adapt the model to a circular geometry, we rewrite \eqref{eq:R2cartesiano} in polar coordinates $(\rho,\theta)$, as follows:
\begin{equation*}
    \begin{cases}
        \displaystyle \dot{\rho} = v_{\rho} \\
        \displaystyle \dot{v}_{\rho} = - \frac{v_\theta^{2}}{\rho} - \frac{\mu}{\rho^2} + \Bar u \cos{\phi}\\
        \displaystyle \dot{\theta} = \frac{v_\theta}{\rho} \\
        \displaystyle \dot{v}_{\theta} = - \frac{v_\rho v_\theta}{\rho} + \Bar u \sin{\phi},
    \end{cases}
\end{equation*}
with the initial conditions
\begin{equation*}
    \begin{cases}
        \rho(0) = \rho_{0} \\
        v_{\rho}(0) = {v_\rho}_{0}\\
        \theta(0) = \theta_{0} \\
        v_{\theta}(0) = {v_\theta}_{0},
    \end{cases}
\end{equation*}
where $\phi \in [0,2\pi)$ is the satellite thruster orientation w.r.t. the radius $\rho$, and $\Bar u\in\{T/m,0\}$.
Now, the running cost \eqref{eq:runcost} is recast in terms of polar coordinates as 
\begin{equation}\label{eq:runcost_polar}
    l(\rho,\theta,v_\rho,v_\theta,\Bar u,\phi) = \alpha \Bar u + \beta (a - \Bar{a})^2 + \gamma \Big((e \cos(\omega) - \Bar{e} \cos(\Bar{\omega}))^2 + (e \sin(\omega) - \Bar{e} \sin(\Bar{\omega}))^2\Big),
\end{equation}
where (see \cite{Pon})
\begin{equation*}
	a = a(\rho,\theta,v_\rho,v_\theta) = -\frac{\mu}{v^2 \rho - 2\mu}
\end{equation*}
\begin{equation*}
	e = e(\rho,\theta,v_\rho,v_\theta) = \sqrt{1 + \frac{2 \mathcal{E}h^2}{\mu^2}}
\end{equation*}
where $\mathcal{E}$ is the energy defined as 
\begin{equation*}
	\mathcal{E} = \mathcal{E}(\rho,\theta,v_\rho,v_\theta) = -\frac{\mu}{\rho} + \frac{v^2}{2}
\end{equation*}
$v$ is the velocity vector, 
\begin{equation*}
	v = \sqrt{v_\rho^2 + v_\theta^2}
\end{equation*}
and $h$ is the angular momentum, $h = \rho v_\theta$.
The argument of perigee is defined by
\begin{equation*}
	\omega = \omega(\rho,\theta,v_\rho,v_\theta) =\theta - \theta_*
\end{equation*}
where $\theta_*$ is the true anomaly,  that we derived as
\begin{equation*}
	\theta_* = \theta_*(\rho,\theta,v_\rho,v_\theta) = 2 \arctan{\bigg(\frac{\rho v_\rho}{\rho(e-1) + p} \sqrt{\frac{p}{\mu}}\bigg)}
\end{equation*}
with $p$ semilatus rectum, $p = a(1 - e^2)$.

\paragraph{Numerical acceleration techniques: Policy Iteration.}
As already remarked, equation \eqref{eq:schema} is a fixed point problem, in which the right-hand side is a contraction with Lipschitz constant $e^{-\lambda\Delta t}$. Therefore, given an initial guess $V^{(0)}$, the {\em value iteration} form
\begin{equation}
v_j^{(k+1)} = \min_{f\in U} \left\{\Delta t \> l(x_j,f) + e^{-\lambda\Delta t} I\left[V^{(k)}\right](x_j+\Delta t\>\Phi(x_j,f))\right\}
\end{equation}
allows to obtain the solution $V$ of \eqref{eq:schema} as $V=\lim_k V^{(k)}$.

However, the convergence of the value iteration can be very slow and, for this reason, we solve \eqref{eq:schema} through the faster \textit{policy iteration}. This technique can be proved to be a Newton-like algorithm, and has superlinear convergence under relatively mild assumptions (in our case, the monotonicity of the scheme) \cite{PutBru,SanRus,BMZ}.

We start by choosing an initial guess $f_j^{(0)}$ for the feedback to be applied at any node $x_j$.  Freezing the control in \eqref{eq:schema},  the Bellman equation becomes linear, and it is solved as an advection equation (this phase is termed as {\em policy evaluation}). Then, a new policy is computed (this step represents the {\em policy improvement}) and a new iteration starts.  Then, the procedure is:
\begin{enumerate}
	\item Set an initial guess $\left\{f_j^{(0)}\right\}_j$ for the policy, and set $k=0$.
	\item For $k\geq 0$, until convergence,  compute for any node of the grid the solution $v_j^{(k)}$ of
	\begin{equation}\label{eq:p_eval}
		v_j^{(k)} = \Delta t \> l\left(x_j,f_j^{(k)}\right) + e^{-\lambda\Delta t} I\left[V^{(k)}\right]\left(x_j+\Delta t\>\Phi\left(x_j,f_j^{(k)}\right)\right)
	\end{equation}
	\item Compute a new policy on the solution $V^{(k)}$ as
	\begin{equation}\label{eq:p_impr}
		f_j^{(k+1)} =\arg \min_f \left\{ \Delta t \> l(x_j,f) + e^{-\lambda\Delta t} I\left[V^{(k)}\right](x_j+\Delta t\>\Phi(x_j,f))\right\}
	\end{equation}	
	\item Increment $k$ and go to 2.
\end{enumerate}
Once the coefficients of the unknowns $v_i^{(k)}$ in the interpolation have been explicitly written (we recall that, for a local monotone interpolation as $\mathbb Q_1$ is, $I[V^{(k)}](x_j+\Delta t\>\Phi(x_j,f))$ is a convex combination of the values $v_i^{(k)}$ in the neighbourhood of $x_j+\Delta t\>\Phi(x_j,f)$), equation \eqref{eq:p_eval} amounts to a linear system,
\[
v_j^{(k)} - e^{-\lambda\Delta t} I\left[V^{(k)}\right]\left(x_j+\Delta t\>\Phi\left(x_j,f_j^{(k)}\right)\right)= \Delta t \> l\left(x_j,f_j^{(k)}\right) \qquad (j=1,\ldots,N).
\]
Given the dimension and the sparsity of the system, the use of efficient sparse solvers like GMRES is essential.

\paragraph{Efficient computation of the minimum.}
A third possible bottleneck of the scheme is in the minimization phase, to be performed in the policy improvement substep \eqref{eq:p_impr}. Since this minimization is repeated for each node of the space grid, speed becomes crucial -- possibly, at the expense of an extreme accuracy. Note that the error in this phase needs not be of higher order than the other terms in the consistency error of the scheme (see \cite{FF94} for an analysis of the error introduced by the inexact computation of the minimum).

The starting point for minimization is a tabulation on the control set, which may be used either for an exhaustive (although not very accurate) search, or to provide a first guess for a more precise optimization algorithm. In the code under consideration, we have used a tabulation with 73 values of the control, i.e., the angles $\phi$ of orientation of the thrust with a $5^\circ$ tabulation step for $\Bar u=T/m$, plus the null thrust $\Bar u=0$.

At the first iteration, we must necessarily minimize over all the possible controls of the tabulation. However, we do not expect that the optimal policy might change much from one iteration to the following. For this reason, at successive iterations we start from the previous optimum value, and walk towards the minimum, one step at a time, until a local minimizer is found. Given the small neighbourhood in which the points $x_j+\Delta t\>\Phi(x_j,f)$ are confined, the function to be minimized is not expected to have an overly complex structure, and detection of non-global minima is unlikely.

Of course, once the minimum point of the tabulation is detected, a more precise computation of the optimal control at $x_j$ can be performed, e.g., by a bisection or golden-section search \cite{Fle}.

\paragraph{Computation of the approximate optimal solution.}
Once the approximate value function has been computed, an approximate optimal feedback can be defined on the basis of the current state $y(t)$ of the (real) system by
\begin{equation}\label{eq:feed}
F_\Delta(y(t)) =\arg \min_f \left\{ \tau \> l(x_j,f) + e^{-\lambda\tau} I[V](x_j+\tau\>\Phi(x_j,f))\right\},
\end{equation}
for a suitable time step $\tau$. This strategy can be proved to provide an asymptotically optimal solution for the real system, as $\Delta t,\Delta x,\tau \to 0$, and $\Delta t=o(\tau)$ (see \cite{FF13}). In the numerical examples, the evolution of the real system will be computed itself in approximate form, via a fourth-order Runge--Kutta method.

\section{Numerical Tests}\label{sec:test}

In this section, we show a couple of numerical tests to validate the proposed numerical scheme. The tests are carried out in two phases: in the first, the value function is computed, while in the second an orbit raising is simulated, with an optimal feedback computed by \eqref{eq:feed}. The parameters are tuned on the case of small satellites with electric propulsion.

\paragraph{Example 1.} In this example, we assume a target orbit having a semi-major axis of 7000 km, an eccentricity of $0.001$ and an argument of the perigee of $0^\circ$. The acceleration given by the maximum thrust has been set to $\bar u=5\cdot 10^{-7}\>\mathrm{km/s^2}$. The parameters of the cost functional are $\alpha = 2.04\cdot 10^{-8}$, $\beta = 2.31\cdot 10^{-2}$, $\gamma = 1.5$ and $\lambda=10^{-3}$. The stopping tolerance of the Policy Iteration is $10^{-7}$.

The computational domain is defined by:
\[
\rho\in [6930,7070]\> \mathrm{km}, \quad \theta\in [0,2\pi), \quad v_\rho\in [-0.01,0.01]\> \mathrm{km/s}, \quad v_\theta\in [7.526,7.566]\> \mathrm{km/s}.
\]
This domain has been discretized with a grid of respectively 110 nodes on $\rho$ and 30 nodes on $\theta$, $v_\rho$ and $v_\theta$ (for a total number of about $3\cdot 10^6$ nodes), and the time step has been set to $\Delta t=10$ s. The computation of the value function on a MacBook Air laptop, for a Mex-compiled Matlab code, requires about 18 minutes.

It is worth comparing the efficiency of the optimized code versus the simpler versions. The CPU time required for a value iteration solver (with the same stopping tolerance) raises to about 15 hours, i.e., about 50 times higher. Performing minimization in \eqref{eq:p_impr} by an exaustive tabulation increases the global CPU time of about 20\%; taking into account that minimization is only a single building block of the scheme, this means anyway a relevant improvement in the efficiency of this phase.

\begin{figure}[t!]
\centering
\includegraphics[width=7cm]{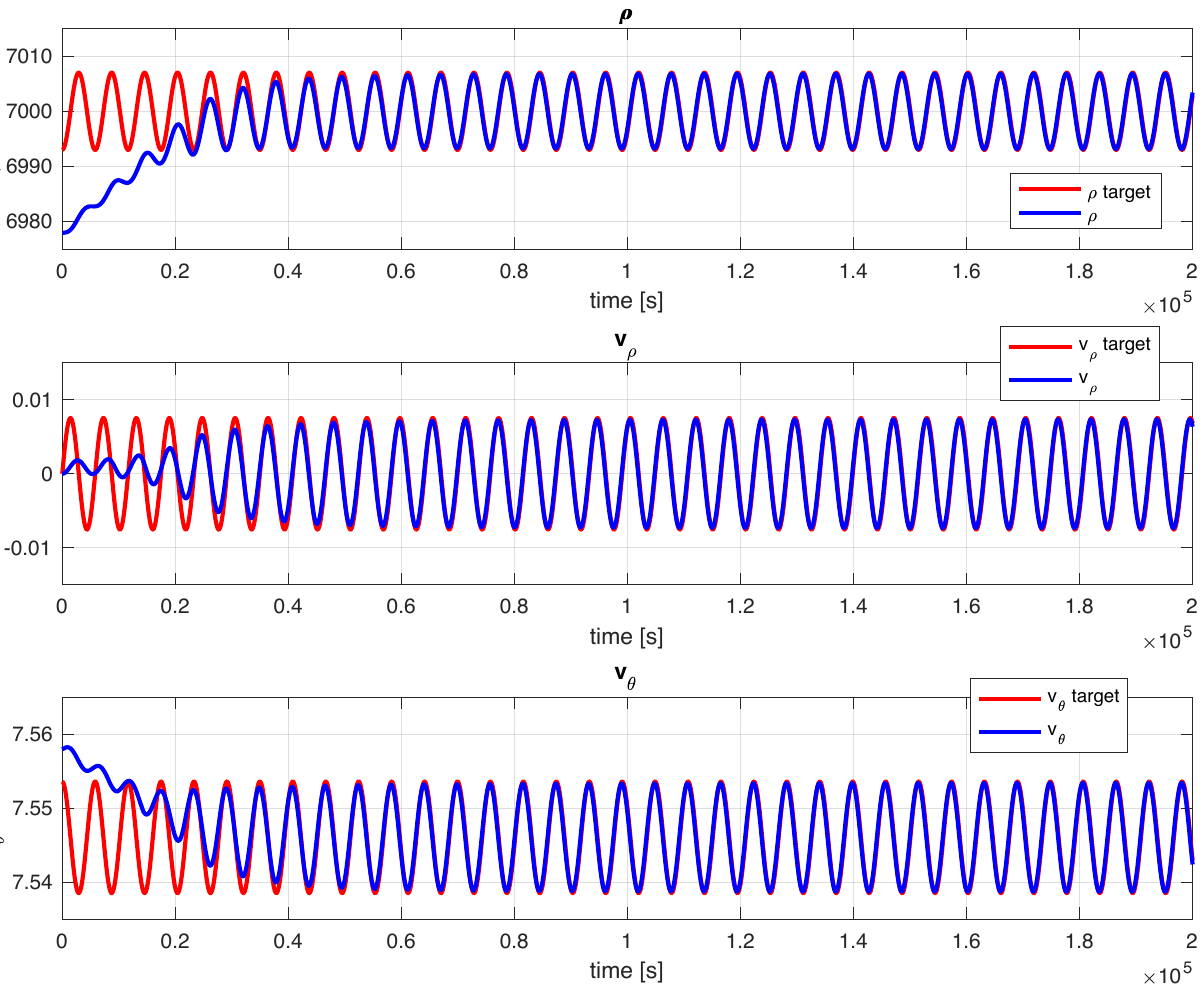} 
\caption{State variables $\rho$, $v_\rho$, $v_\theta$ for the approximate optimal trajectory, compared to the target state variables}\label{fig:f2}
\end{figure}

\begin{figure}[t!]
\centering
\includegraphics[width=7cm]{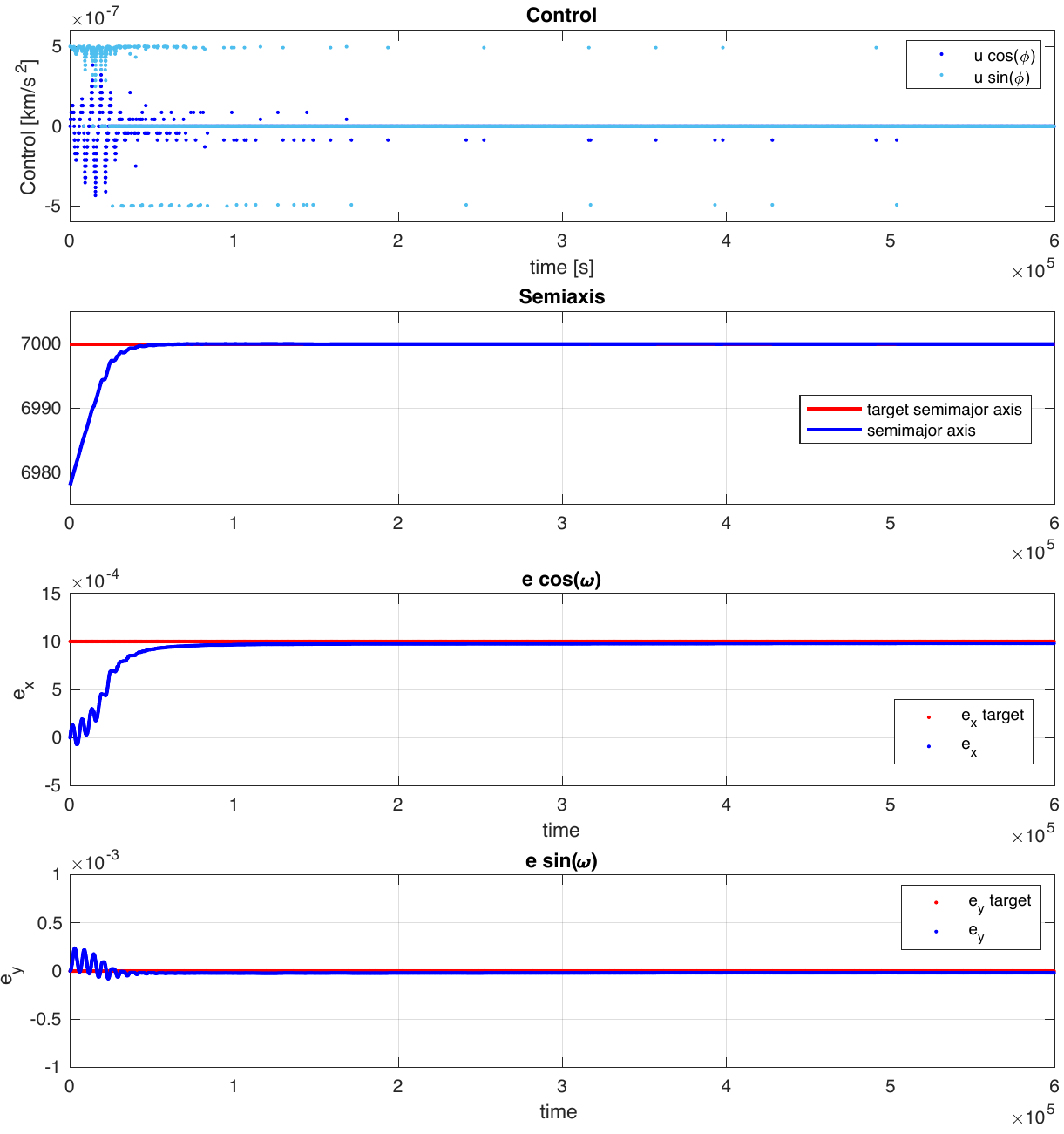}
\caption{Control and orbital parameters $a$, $e_x$, $e_y$ for the approximate optimal trajectory, compared to the target orbital parameters}\label{fig:f3}
\end{figure}

\begin{figure}[t!]
\centering
\includegraphics[width=7cm]{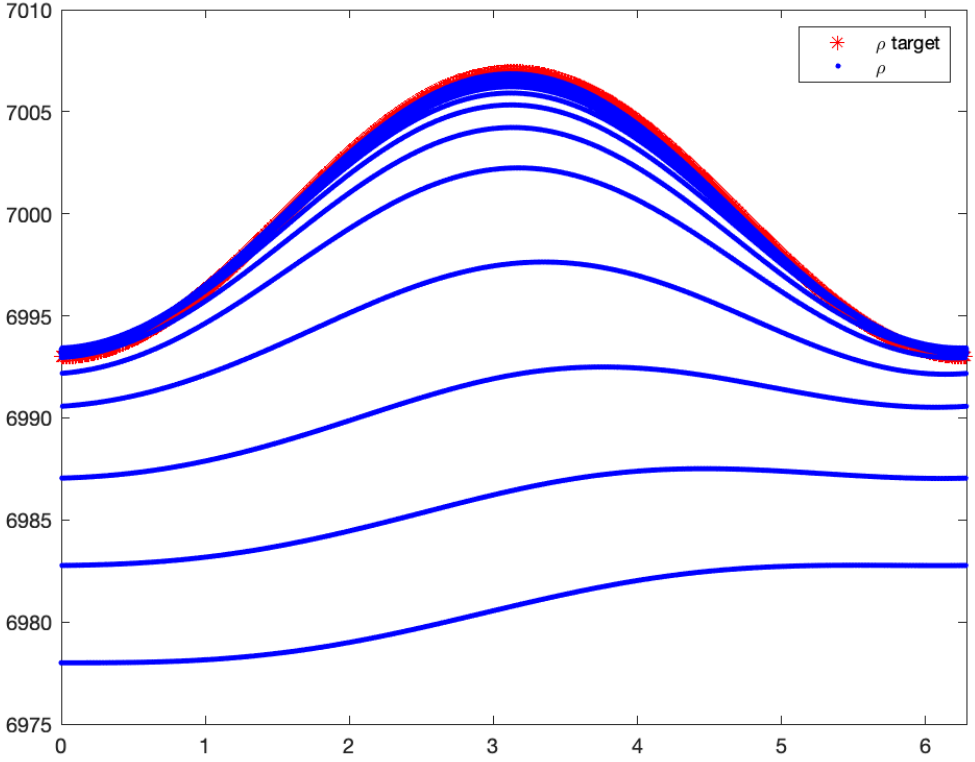}
\includegraphics[width=7cm]{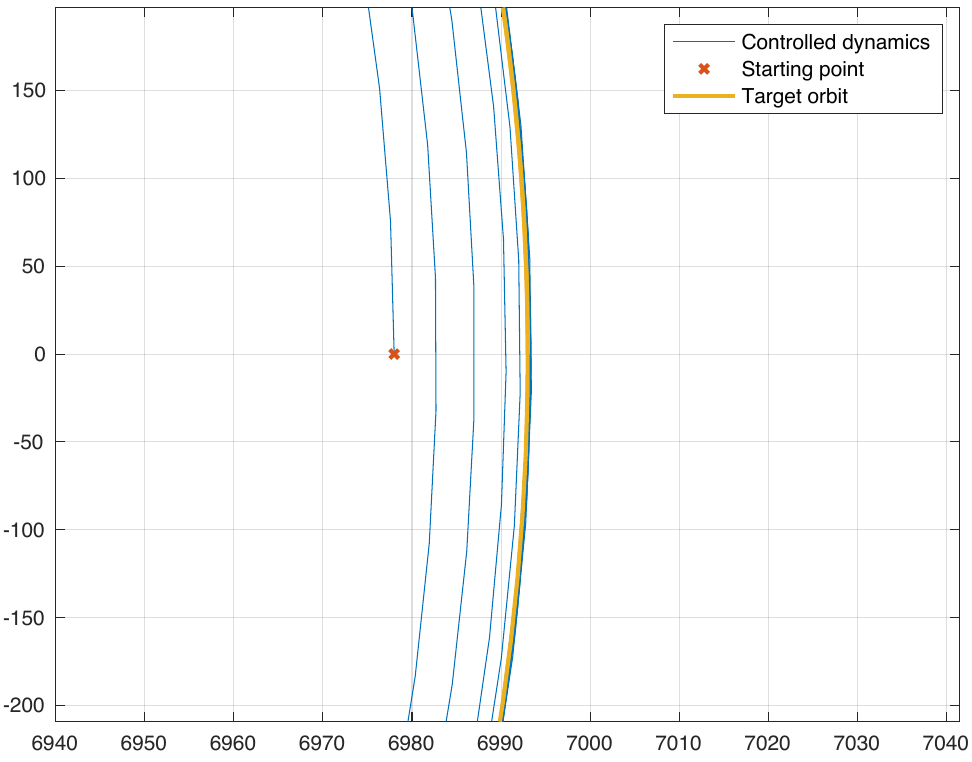}
\caption{State variable $\rho$ as a function of $\theta$ for the approximate optimal trajectory, compared to the $\rho(\theta)$ of the target orbit (left) and a zoom of the approximate optimal trajectory in cartesian coordinates (right)}\label{fig:f4}
\end{figure}

After the computation of value function and optimal feedback, we have simulated an orbit raising manoeuvre starting from a circular preliminary orbit with $\rho\equiv 6978$ km towards the target orbit above. We show in Fig. \ref{fig:f2} the evolution of the state variables  $\rho$, $v_\rho$ and $v_\theta$ versus the target state variables, whereas Fig. \ref{fig:f3} performs the same comparison for the orbital parameters  $a$, $e_x$ and $e_y$. Last, the left plot in Fig. \ref{fig:f4} we show the evolution of $\rho$ as a function of $\theta$ during the orbit raising, compared with the $\rho/\theta$ relationship of the target orbit, while the right plot provides a zoom of the approximate optimal orbit in a neighbourhood of the initial point. The plots of the approximate optimal solution show the expected behaviour: the satellite raises from the preliminary orbit to the target orbit, then switches the control off. The state reaches the neighbourhood of the target orbit within a time of about $10^5$ s (i.e., about 28 hours) and, after further adjustments, the thruster is switched off after about $5\cdot 10^5$ s (i.e., about 6 days), with an error of the order of $10^{-2}$ km on the radius of the orbit. 


\paragraph{Example 2.} In this second example, we test the robustness of the approximate optimal feedback with respect to perturbations and missing terms in the model of the physical system. In practice, we compute the value function and the optimal feedback on the basis of the simplified model \eqref{eq:R2cartesiano}, but apply this feedback to a physical system in which some of the additional perturbations in \eqref{eq:model0} appear: in particular, we introduce in the ``real'' system the aerodynamic drag term.
This term, which should in principle depend on various parameters, among which the orientation of the orbit plane with respect to the Earth and the direction of rotation, is split into tangential and radial components, and approximated as
\[
{a_D}_\rho = - C v_\rho^2, \quad {a_D}_\theta = - C v_\theta^2,
\]
with
\[
C = \frac{C_D S}{2m} \rho_D,
\]
in which we have assumed that the relative speed of the satellite with respect to the atmosphere coincides with the satellite speed itself, and use a constant density $\rho_D=8.09\cdot 10^{-14}$ kg/m$^3$
for the air, given the small altitude interval. The other parameters are set as $C_D=2.2$, $S=2.25$ m$^2$ and $m=350$ kg.

\begin{figure}[t!]
\centering
\includegraphics[width=12cm]{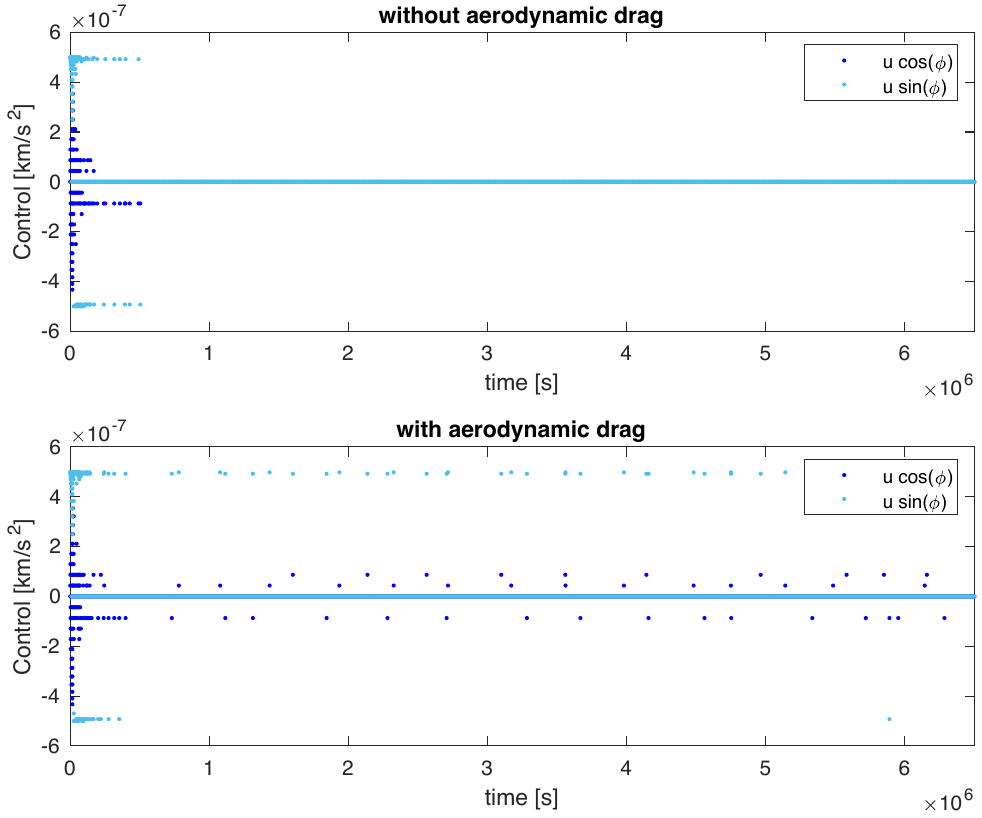}
\caption{Control $u$ as a function of $t$ for the approximate optimal trajectory, without aerodynamic drag (upper) and in presence of it (lower)}\label{fig:f5}
\end{figure}

\begin{figure}[t!]
\centering
\includegraphics[width=12cm]{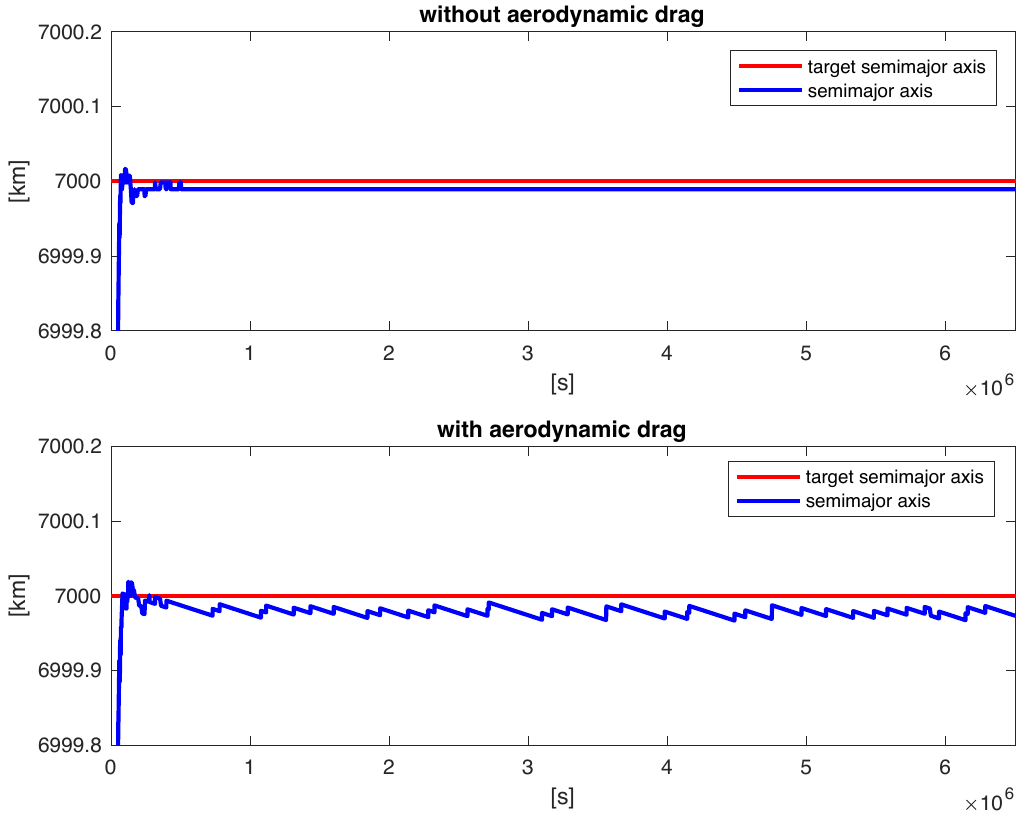}
\caption{Semi-major axis $a$ as a function of $t$ for the approximate optimal trajectory, without aerodynamic drag (upper) and in presence of it (lower)}\label{fig:f6}
\end{figure}

The approximate optimal feedback provides a solution very close to that of the previous example, except for some adjustment manoeuvres. These are shown in Figures \ref{fig:f5}--\ref{fig:f6}, in which we compare the approximate optimal control of the system without drag (Fig. \ref{fig:f5}, upper plot) with the corresponding control in presence of drag (Fig. \ref{fig:f5}, lower plot), and the semi-major axes in the two cases (Fig. \ref{fig:f6}). Here, the state reaches the target orbit and switches off the control, but, for the system with aerodynamic drag, new adjustments of the trajectory are performed at later times to correct the apparent energy dissipation. 

\section*{Conclusions}

We have proposed, and preliminary validated through numerical testing, a DP-based approach for the optimal control in orbit transfer problems, with a special focus on LEO satellites. By choosing an infinite horizon cost functional, we obtain a static feedback control, which may work as a stabilizing feedback around the target orbit, and proves to be robust enough with respect to perturbations of the state, including the presence of terms neglected in the mathematical model. In addition, we can compute a fairly accurate solution despite the relatively coarse discretization. In a forthcoming step, solutions computed via DP will be compared with that resulting from the application of open-loop techniques, for a more quantitative assessment.

In this preliminary study, we have tested in a planar situation various techniques that allow for relevant reductions of the computational effort (which remains a major drawback of DP methods); a clever use of these tricks may prepare for efficiently tackling the full $\mathbb R^6$ problem of orbit transfer.

\section{Acknowledgements}

R.F. and G.V. are members of INdAM--GNCS, and have been partially supported by the PRIN project "Optimal control problems: analysis, approximation and applications".


\begin{thebibliography}{20}

\bibitem{BCD} M. Bardi, I. Capuzzo Dolcetta, Optimal control and viscosity Solutions of Hamilton--Jacobi--Bellman equations, Birkh\"auser, Boston, 1997.

\bibitem{BS} G. Barles, P.E. Souganidis, {\em Convergence of approximation schemes for fully nonlinear second-order equations}, Asymptotic Analysis, {\bf 4} (1991), 271--283.  

\bibitem{B} R. Bellman, Dynamic Programming, Princeton University Press, Princeton NJ, 1957.

\bibitem{BBDZ21} O. Bokanowski, E. Bourgeois, A. D\'esilles, H. Zidani, {\em Global optimization approach for the ascent problem of multi-stage launchers}, in: H.G. Brock, W. J\"ager, E. Kostina, H.X. Phu (eds.) Modeling, Simulation and Optimization of Complex Processes HPSC 2018, Springer, Cham, 2021, 1--42.

\bibitem{BMZ} O. Bokanowski, S. Maroso, H. Zidani, {\em Some convergence results for Howard's algorithm}, SIAM Journal on Numerical Analysis {\bf 47} (2009), 3001--3026.

\bibitem{BBDZ15} E. Bourgeois, O. Bokanowski, A. D\'esilles, H. Zidani, {\em Optimization of the launcher ascent trajectory leading to the global optimum without any initialization: the breakthrough of the HJB approach}, 6th European Conference for Aeronautics and Space Sciences (EUCASS), Vol. 29, 2015.

\bibitem{BBZD} E. Bourgeois, O. Bokanowski, H. Zidani, A. D\'esilles, {\em New improvements in the optimization of the launcher ascent trajectory through the HJB approach}, 7th European Conference for Aeronautics and Space Sciences (EUCASS), 2017.

\bibitem{CFTZ} J.-B. Caillau, R. Ferretti, E. Trélat, H. Zidani, {\em An algorithmic guide for finite-dimensional optimal control problems}, Handbook of Numerical Analysis, {\bf 24} (2023), 559--626.

\bibitem{CriMar} E. Cristiani, P. Martinon, {\em Initialization of the shooting method via the Hamilton--Jacobi--Bellman approach}, J. Optim. Theory Appl. {\bf 146} (2010), 321--346.

\bibitem{Cur} H.D. Curtis, Orbital mechanics for engineering students: revised reprint, Butterworth--Heinemann, 2020.

\bibitem{FF94} M. Falcone, R. Ferretti, {\em Discrete time high-order schemes for viscosity solutions of Hamilton--Jacobi--Bellman equations}, Num. Math. {\bf 67} (1994), 315--344.

\bibitem{FF13} M. Falcone, R. Ferretti, Semi-Lagrangian approximation schemes for linear and Hamilton--Jacobi equations, SIAM, Philadelphia, 2013.

\bibitem{Fle} R. Fletcher, Practical methods of optimization, John Wiley \& Sons, 2000.

\bibitem{Pon} M. Pontani, Spaceflight Mechanics lecture notes, Sapienza Universit\`a di Roma, Rome, Italy, 2018.

\bibitem{PutBru} M.L. Puterman, S.L. Brumelle, {\em On the convergence of policy iteration in stationary dynamic programming}, Mathematics of Operational Research {\bf 4} (1979), 60--69.

\bibitem{SanRus} M.S. Santos, J. Rust, {\em Convergence properties of policy iteration}, SIAM J. on Control and Optimization {\bf 42} (2004), 2094--2115.

\end{thebibliography}
\end{document}